\documentclass[12pt]{elsarticle}
\makeatletter
\def\ps@pprintTitle{%
 \let\@oddhead\@empty
 \let\@evenhead\@empty
 \def\@oddfoot{}%
 \let\@evenfoot\@oddfoot}
\makeatother

\usepackage{amssymb}
\usepackage{amsmath}
\usepackage{tabls}

\newtheorem{thm}{Theorem}

\newtheorem{cor}{Corollary}
\newtheorem{prp}{Proposition}

\newtheorem{obs}{Observation}

\newtheorem{prt}{Property}

\newtheorem{exmp}{Example}

\newproof{pf}{Proof}

\journal{}

\begin{document}

\begin{frontmatter}



\title{Some properties of $\{k\}$-packing function problem in graphs}


\author[frs]{Jozef J. Kratica}
\ead{jkratica@mi.sanu.ac.rs}

\author[sec]{Aleksandar Lj. Savi\'{c}}
\ead{aleks3rd@gmail.com}

\author[thd]{Zoran Lj. Maksimovi\'{c}\corref{cor1}}
\ead{zoran.maksimovic@gmail.com}

\address[frs]{
Mathematical Institute,
Serbian Academy of Sciences and Arts, 
Kneza Mihaila 36/III, 11000 Belgrade, Serbia}

\address[sec]{
Faculty of Mathematics, University of Belgrade, Studentski
trg 16/IV, 11000 Belgrade, Serbia}

\address[thd]{
University of Defence, Military Academy, Generala Pavla
Juri\v{s}i\'{c}a \v{S}turma 33, 
11000 Belgrade, Serbia}
\cortext[cor1]{Corresponding author}

\begin{abstract}
	The recently introduced  $\{k\}$-packing function problem is
	considered in this paper. Special relation between a case when $k=1$, $k\ge 2$ and linear programming
	relaxation is introduced with sufficient conditions for optimality.
	For arbitrary simple connected graph $G$ there is construction procedure 
	for finding values of $k$ for which $L_{\{k\}}(G)$ can be determined in the polynomial time. 
	Additionally, relationship between 
	$\{1\}$-packing function and independent set number is established.
	Optimal values for some special classes of graphs
	and general upper and lower bounds are introduced.
\end{abstract}


\end{frontmatter}



\section{Introduction}

In this paper, we will consider simple, finite and undirected graphs. 
For a graph $G$, $V(G)$ and $E(G)$ denote its vertex and edge
sets, respectively. Further, for any $v \in V(G)$, 
its open neighborhood $N_G(v)$ 
is the set  of all vertices that are adjacent
to $v$, and its closed neighborhood is $N_G[v]=N_G(v)\cup \{v\}$. 
For a function $f:V(G)\rightarrow {\mathbb N}\cup\{0\}$, and $A\subseteq V(G)$
it should be denoted $f(A)=\sum\limits_{v\in A}f(v)$.
Let $|V(G)|=n$ and $A_G=[a_{ij}]_{n\times n}$ where 
$a_{ij}=
\begin{cases}
	1, i=j \lor (i,j)\in E(G) \\
	0, { otherwise}           
\end{cases}$

For a graph $G$ and a positive integer $k$, a function $f:V(G)\rightarrow {\mathbb N}\cup\{0\}$, is a $\{k\}$-packing function of graph $G$,
if for each vertex $v\in V(G)$ value $f(N_G[v])$ is at most $k$.
The maximum possible value of $f(V(G))$ over all $\{k\}$-packing functions of graph $G$
is denoted as $L_{\{k\}}(G)$.
Formally, $L_{\{k\}}(G)=\max\limits_{f:V(G)\rightarrow {\mathbb N}\cup\{0\}} \{f(V(G))|(\forall v\in V(G))f(N_G[v])\le k\}$.

The distance between vertices $u$ and $v$, denoted as $d_G(u,v)$ is the length of the shortest $u-v$
path. The square of a graph $G$, named $G^2$, is the graph obtained from $G$ by adding all edges between vertices 
from $V(G)$ that have a common neighbor, i.e. $G^2 = (V(G),E(G^2))$, where $E(G^2) = \{(u,v) \in V(G) \times V(G)\,|\, d_G(u,v) \le 2\}$.
The complement of a graph $G$, named $\overline{G}$, is defined as  $\overline{G} = (V(G), \overline{E(G)})$, where
$\overline{E(G)} = \{(u,v) \in V(G) \times V(G)\,|\, u \ne v \,\wedge\, (u,v) \notin E(G)\}$. 
The independent set $I(G)$ of a graph is a set of vertices, subset of $V(G)$, such that there are no edges between them,
i.e. ($u,v \in I(G) \Rightarrow (u,v) \notin E(G)$). Independence number of a graph,
named $\alpha(G)$ is the cardinality of a maximal independent set $I(G)$.

For $k$ being fixed positive integer Meir and Moon \cite{meir75} introduced $k$-packing set $P\subset V(G)$ as 
a set of vertices such that distance between $u$ and $v$ is greater than $k$ 
for distinct $u,v\in P$, and 
$k$-packing number ($\rho_k(G)$) as the number of vertices of such largest set. 
It stands that $\rho_1(G)=\alpha(G)$ is the independence number.

Gallant et al. in \cite{gal10} introduced $k$-limited packing as a modification of packing number problem allowing
that intersection of each closed neighborhood with a given set contains no more than $k$ vertices. 
In \cite{dob10, dob11} Dobson et al. proved that $k$-limited packing is NP-complete for
split and bipartite graphs. It was also shown that $P_4$-tidy graphs are solvable in 
polynomial time.

The notion of $\{k\}$-packing function was introduced by Leoni and Hinrichsen \cite{leo14} 
as a variation of $k$-limited packing in order to solve the problem of locating garbage dumps
in a given city. In this scenario, it is possible to place more than one dump in a certain location, 
requesting that no more than $k$ dumps are placed in each vertex and its neighborhood. 
Although notation is similar, for $k\ge 2$ it must be clearly distinguished $k$-limited packing function and $\{k\}$-packing function. 
Relationship between $k$-limited packing and $\{k\}$-packing function 
is established 
in \cite{leo16}. It was stated that $L_{\{k\}}(G)\ge L_k(G)$.
Additionally,
in \cite{leo14} where it is shown that $L_{\{k\}}(G)=L_k(G\otimes K_k)$
($\otimes$ is a strong product of graphs). 

\begin{prp} \label{a1}
	{\rm (\cite{dob17})}
	For a graph $G$ and a positive integer $k$ it holds 
	$L_{\{k\}}(G)\ge k\cdot L_1(G)$	
\end{prp}

\begin{prp} \label{m1}
	{\rm (\cite{moj17})}
For any connected graph
$G$ and integer $k\in \{1,2\}$
$L_{k}(G)\ge \lceil \frac{k\cdot diam(G)+k}{3}\rceil$
\end{prp}

\begin{prp}
	\label{path}
	{\rm (\cite{gal10,fren09,dob17})} For path $P_n$ holds $L_{\{k\}}(P_n)=\lceil\frac{n}{3}\rceil\cdot k$.
\end{prp}
\begin{pf} The proposition directly holds from the following statements: 	\\
	\begin{itemize}
		\item In (\cite{gal10}) in Lemma 3 it was proven that $L_1(P_n)=\lceil\frac{n}{3}\rceil$;
		\item From \cite{fren09} Theorem 1. it holds that 
		$\gamma(P_n)=\lceil\frac{n}{3}\rceil$;
		\item Finally in (\cite{dob17}) Theorem 3.1
		it was proven that  
		$\gamma(G)=L_1(G)\Rightarrow L_{\{k\}}(P_n)=k\cdot L_1(P_n)$.
	\end{itemize}

	Therefore
	$L_{\{k\}}(P_n)=k\cdot L_1(P_n)=\lceil\frac{n}{3}\rceil\cdot k$. 
\end{pf}

\begin{thm}\label{np}\mbox{\rm \cite{dob17}}
	The $\{k\}$-packing function problem is NP-complete for all integer $k$ fixed.
\end{thm}

The polynomial equivalence between $\{k\}$-packing function problem and 
$k$-limited packing in graphs is discussed in \cite{leo16}.


\section{New theoretical properties}

In this section, relationship between $\{k\}$-packing,
$\{1\}$-packing problem 
and relaxation of $\{1\}$-packing will be established as well as 
some properties of $\{k\}$-packing function problem for certain 
classes of graphs. Without loss of generality, we will assume
that considered graphs are connected and have at least two vertices since 
if the graph is not connected we can consider connected components instead,
using the following simple property.

\begin{prt}
	If $G$ is not connected and has connected components $Con_1(G)$, $Con_2(G)$, \ldots $Con_{nc}(G)$ then
	$L_{\{k\}}(G)=\sum\limits_{i=1}^{nc} L_{\{k\}}(Con_i(G))$
\end{prt}

\begin{pf}
	Let $v\in V$ be an arbitrary vertex from a connected component $Con_j(G)$.
	Since $v\in Con_j(G)\Rightarrow N[v]\subseteq Con_j(G)$, then all constraints
	$f(N_G[v])\le k$ can be grouped by connected components and considered independently.
\end{pf}

Let $Z^*_{rlx}(G)$ be an optimal solution of the relaxed
$\{1\}$-packing problem. 
Relaxation is performed by 
$Z^*_{rlx}(G)=\max\limits_{f:V(G)\rightarrow [0,+\infty)} \{f(V(G))|(\forall v\in V(G))f(N_G[v])\le 1\}$,
	i.e. relaxed packing function can take 
	fractional (real) values.  
	
	Now we can formulate simple, but effective, relation among $L_{\{k\}}(G)$, $L_{\{1\}}(G)$ and $Z^*_{rlx}(G)$.
	
	\begin{prp} 
		\label{prop1}For arbitrary $k\in \mathbb{N}$ it stands 
		$ k\cdot L_{\{1\}}(G) \le L_{\{k\}}(G) \le k\cdot Z^*_{rlx}(G)$
	\end{prp}
	\begin{pf}
		It should be noted that proof cannot be based on Proposition \ref{a1} and
		fact that 
		$L_{\{k\}}(G)\ge L_k(G)$. 
		
		Let $f$ be a $\{1\}$-packing function of $G$ with the maximum value of all such functions. Then function $g:V(G)\rightarrow \mathbb{N}\cup \{0\}$ 
		such that $g(v)=k\cdot f(v)$ is obviously a $\{k\}$-packing function of $G$. 
		Consequently, 
		$L_{\{k\}}(G)=\max\limits_{h:V(G)\rightarrow {\mathbb N}\cup\{0\}} \{h(V(G))|(\forall v\in V(G))h(N_G[v])\le k\}
		\ge g(V(G))$. Therefore, 
		$L_{\{k\}}(G) \ge k\cdot L_{\{1\}}(G)$.
		
		It should be noted that $L_{\{k\}}(G) \ge k\cdot L_{\{1\}}(G)$
		directly follows from Proposition \ref{a1} and
		fact that $L_{\{1\}}(G)= L_1(G)$. 
	
		Let $f_{rlx}:V(G)\rightarrow [0,+\infty)$ be a relaxed $\{1\}$-packing function with maximum value of all such functions. 
			As it stands that
			$$(\forall v\in V(G))f_{rlx}(N_G[v])\le 1 \Rightarrow k\cdot f_{rlx}(N_G[v])\le k$$
			
			\noindent
			and $\{k\}$-packing function has non negative integer values,
			then 
			
			$L_{\{k\}}(G)=\max\limits_{h:V(G)\rightarrow {\mathbb N}\cup\{0\}} \{h(V(G))|(\forall v\in V(G))h(N_G[v])\le k\} \le$
			
			$\le k\cdot \max\limits_{f:V(G)\rightarrow [0,+\infty)} \{f(V(G))|(\forall v\in V(G))f(N_G[v])\le 1\}=k\cdot Z^*_{rlx}(G)$.
				\end{pf}			
	
				It is interesting to find when equalities hold, i.e. when $k\cdot L_{\{1\}}(G) = L_{\{k\}}(G)$ or $k\cdot L_{\{1\}}(G)= k\cdot Z^*_{rlx}(G)$. 
				Sufficient condition for both equalities will be given in the following theorem.
				
				\begin{thm}
					\label{main}
					If $A_G$ is a totally unimodular matrix, then $L_{\{k\}}(G) = k\cdot L_{\{1\}}(G)= k\cdot Z^*_r(G)$ holds.
				\end{thm}
				
				\begin{pf}
					Let $G=(V,E)$ be a graph whose $A_G$ is a totally unimodular matrix. Let us consider $\{k\}$-packing function problem. 
					The problem can be formulated as a following integer linear program. 
					Let us denote the variables $x_i, \, i=1, \ldots, |V|$ 
					such that $x_i=f(i)$. Then, $\{k\}$-packing function problem can be formulated as 
					
					\begin{equation}
						\max \sum\limits_{i = 1}^{|V|} {x_i } 
					\end{equation}
					
					subject to
					
					\begin{equation}
						\sum\limits_{j \in N_G [i]} {x_j }  \le k,\quad i = 1, \ldots ,|V| 
					\end{equation}
					
					\begin{equation}
						x_i  \in \{ 0, 1, \ldots ,k\} ,\quad i = 1, \ldots ,|V|
					\end{equation}

					It is easy to see that condition $\sum\limits_{i \in N_G [j]} {x_i }  \le k$ 
					could be replaced with $\sum\limits_{j = 1}^{|V|} {a_{ij} x_j }  \le k$ where $a_{ij}$ are elements of matrix $A_G$. 
					Now, the formulation is

					\begin{equation}
						\label{ILPobj}
						\max \sum\limits_{i = 1}^{|V|} {x_i } 
					\end{equation}
					
					subject to
					
					\begin{equation}
						\sum\limits_{j = 1}^{|V|} {a_{ij} x_j }  \le k,\quad i = 1, \ldots ,|V|
					\end{equation}
					
					\begin{equation}
						\label{ILPvar}
						x_i  \in \{0, 1, \ldots ,k\} ,\quad i = 1, \ldots ,|V|
					\end{equation}

					
					Since this is Integer Linear Programming (ILP) formulation, it is natural to consider its relaxation. 
					Instead of integer constraint $x_i  \in \{ 1, \ldots ,k\}$, let us consider non-negativity constraint $x_i \ge 0$. 
					From the first constraint, it is obvious that for every vertex $i$ will be $x_i \le k$. Let us now consider linear programming formulation 
					
					\begin{equation}
						\label{LPobj}
						\max \sum\limits_{i = 1}^{|V|} {x_i } 
					\end{equation}
					
					subject to

					\begin{equation}
						\sum\limits_{j = 1}^{|V|} {a_{ij} x_j }  \le k,\quad i = 1, \ldots ,|V|
					\end{equation}
					
					\begin{equation}
						\label{LPvar}
						x_i \ge 0  ,\quad i = 1, \ldots ,|V|
					\end{equation}

					
					Note that this formulation for $k=1$ is exactly Linear Programming (LP) formulation of
					$Z^*_{rlx}(G)$: \\

					\begin{equation}
						\label{LPrelobj}
						\max \sum\limits_{i = 1}^{|V|} {x_i } 
					\end{equation}
					
					subject to

					\begin{equation}
						\sum\limits_{j = 1}^{|V|} {a_{ij} x_j }  \le 1,\quad i = 1, \ldots ,|V|
					\end{equation}
					
					\begin{equation}
						\label{LPrelvar}
						x_i \ge 0  ,\quad i = 1, \ldots ,|V|
					\end{equation}

					Since at least one feasible solution of the formulation above exists, $x_i=0, \, i=1,\ldots, |V|$, and all variables have upper bound, an optimal solution also exists. 
					From the theory of integer linear programming, it is known that polyhedron $X(b)$, defined as $X(b)=\{x|Ax \ge b\}$ for any integer vector $b$, 
					is an integer if and only if the matrix $A$ is totally unimodular. 
					Since polyhedron of relaxation of our problem is $X(b)=\{x| A_Gx \le k \cdot e_{|V|}\}$, 
					where $e_{|V|}=(1, \ldots, 1)^T$ is vector of ones and dimension equal to $|V|$, 
					has totally unimodular matrix $A_G$, it can be concluded that all of polyhedron nodes are integer. 
					This means that all optimal solutions of the relaxation problem are integer. 
					As ILP and LP formulations differ only in the condition of integrality, it can be concluded that 
					optimal solutions of the relaxation and ILP formulation are the same
					under the conditions of this theorem.
					
					We have proved that $L_{\{1\}}(G)=Z^*_{rlx}(G)$. From Proposition \ref{prop1} which
					states that 
					$k\cdot L_{\{1\}}(G) \le L_{\{k\}}(G) \le k\cdot Z^*_{rlx}(G)$ 
					and equality of the first and the third term 
					directly holds $k\cdot L_{\{1\}}(G) = L_{\{k\}}(G) = k\cdot Z^*_{rlx}(G)$.
				\end{pf}

				From the well-known fact that any LP problem has a polynomial complexity, the following assertion holds.
				
				\begin{cor}
					If $A_G$ is a totally unimodular matrix, then $\{k\}$-packing function problem can be solved in polynomial time.
				\end{cor}

				However, total unimodularity of matrix $A_G$ is not necessary condition for 
				$k\cdot L_{\{1\}}(G) = L_{\{k\}}(G) = k\cdot Z^*_{rlx}(G)$ to hold, which is illustrated by the following example.
				
				\begin{exmp}
					Let graph $G$ be a claw graph with four vertices, i.e. $G=(V,E)$,
					where $V=\{1,2,3,4\}$ and $E=\{\{1,2\},\{1,3\},\{1,4\}\}$. 
					Matrix $A_G$ is not totally unimodular 
					since $det(A_G)=-2$.
					Since $N[1]=V(G)$ taking into consideration 
					$L_{\{1\}}(G)$ we have $f(V(G))=f(N[1])\le 1$.
					We can construct $\{1\}$-packing function $f$ where
					$f(V(G))=1$: $f(1)=1$ and $f(2)=f(3)=f(4)=0$.
					It is obvious that constructed function $f$ is also
					maximum $Z^*_{rlx}(G)$ of the relaxation problem.
					From the previous facts, clearly 
					$L_{\{1\}}(G)=Z^*_{rlx}(G)=1$.
					Therefore, by Proposition \ref{prop1} it holds
					$k\cdot L_{\{1\}}(G) = L_{\{k\}}(G)= k\cdot Z^*_{rlx}(G)=k$.
				\end{exmp}

				The following example illustrates the case when 
				$k\cdot L_{\{1\}}(G) < L_{\{k\}}(G)$.
				
\begin{figure}[h]
					\centering\setlength\unitlength{1mm}
					
					\setlength\unitlength{1mm}
					\includegraphics[width=6cm]{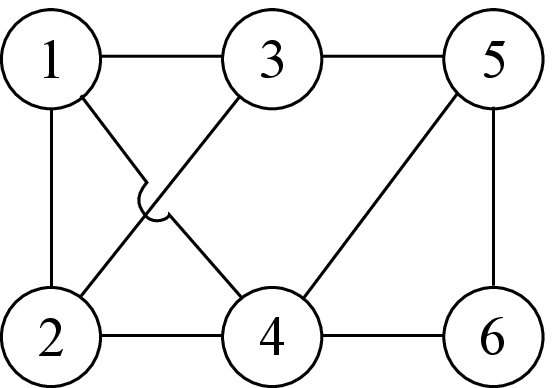}%
					\caption{An example of a graph $G$ where 
						$k\cdot L_{\{1\}}(G) < L_{\{k\}}(G)$.}{}
					\label{f1}
				\end{figure}

				\begin{exmp}
					\label{ex1}Let us consider graph $G$ given in Figure \ref{f1}.
				\end{exmp}
				
				For graph $G$ presented in Figure \ref{f1},
				$2\cdot L_{\{1\}}(G)=2<L_{\{2\}}(G)=3$ holds, 
				since values $ L_{\{1\}}(G)=1$ and $ L_{\{2\}}(G)=3$ are obtained by a total enumeration.
				For $k=1$,  $\{1\}$-packing function with maximal value is defined as follows: $f_1(1)=1$; $\,f_1(2)=f_1(3)=f_1(4)=f_1(5)=f_1(6)=0$.
				For $k=2$,  $\{2\}$-packing function with maximal value is defined as follows: $f_2(2)=f_2(3)=f_2(6)=1$; $\,f_2(1)=f_2(4)=f_2(5)=0$.

				Next, it will be presented an example where  
				$L_{\{k\}}(G) < k\cdot Z^*_R(G)$.
				
				\begin{exmp} Let graph $G$ be given  $V(G)=\{1,2,\ldots,30\}$			
					and adjacency matrix $A_G$ given in Figure \ref{ag}.
				For graph $G$ presented in Figure \ref{ag},
				$L_{\{1\}}(G)= 1 < \lfloor Z^*_{rlx}(G)\rfloor=2$
				holds. 
				Values $ L_{\{1\}}(G)=1$ can be obtained
					using ILP formulation (4)-(6),
					while 
					$Z^*_{rlx}(G)=\frac{7}{3}$
					can be obtained from relaxed LP formulation (10)-(12).
					Values of function $f$ which correspond
					to $Z^*_{rlx}(G)$
					are: 
					$f(3)=f(4)=\frac{1}{7}$; $\,f(5)=f(6)=\frac{2}{21}$; 
					$\,f(8)=f(13)=f(20)=f(24)=\frac{4}{21}$;
					$\,f(10)=f(16)=f(18)=\frac{5}{21}$;
					$\,f(19)=\frac{1}{21}$ and $f(11)=\frac{1}{3}$.
					For any other vertex $v$, $f(v)=0$.
				\end{exmp}
			
				\begin{figure}[h]
				\centering\setlength\unitlength{1mm}
				\tiny
				\setlength\unitlength{1mm}
						\setcounter{MaxMatrixCols}{30}
$A_G=\begin{bmatrix}
0 & 0 & 0 & 0 & 1 & 0 & 0 & 0 & 0 & 0 & 0 & 0 & 0 & 0 & 0 & 0 & 0 & 1 & 1 & 1 & 1 & 1 & 1 & 0 & 0 & 0 & 1 & 1 & 1 & 0\\
0 & 0 & 0 & 1 & 1 & 1 & 0 & 0 & 0 & 0 & 0 & 1 & 0 & 0 & 1 & 1 & 0 & 1 & 0 & 0 & 1 & 0 & 1 & 1 & 1 & 0 & 0 & 0 & 0 & 0\\
0 & 0 & 0 & 0 & 1 & 0 & 0 & 1 & 0 & 0 & 0 & 1 & 0 & 0 & 0 & 0 & 1 & 1 & 0 & 0 & 0 & 1 & 0 & 1 & 0 & 1 & 0 & 0 & 1 & 0\\
0 & 1 & 0 & 0 & 0 & 0 & 1 & 0 & 0 & 0 & 0 & 0 & 0 & 0 & 1 & 0 & 0 & 0 & 1 & 0 & 1 & 0 & 0 & 0 & 0 & 0 & 0 & 0 & 1 & 0\\
1 & 1 & 1 & 0 & 0 & 0 & 1 & 0 & 1 & 0 & 0 & 0 & 1 & 1 & 1 & 0 & 1 & 1 & 0 & 0 & 0 & 1 & 1 & 1 & 0 & 0 & 0 & 1 & 0 & 1\\
0 & 1 & 0 & 0 & 0 & 0 & 0 & 1 & 1 & 0 & 1 & 0 & 1 & 1 & 0 & 0 & 0 & 0 & 0 & 1 & 0 & 0 & 1 & 0 & 1 & 0 & 0 & 0 & 0 & 1\\
0 & 0 & 0 & 1 & 1 & 0 & 0 & 0 & 1 & 0 & 1 & 0 & 1 & 0 & 1 & 0 & 1 & 0 & 1 & 0 & 0 & 0 & 0 & 1 & 0 & 1 & 0 & 0 & 1 & 1\\
0 & 0 & 1 & 0 & 0 & 1 & 0 & 0 & 0 & 0 & 0 & 0 & 0 & 1 & 0 & 0 & 1 & 0 & 0 & 0 & 0 & 0 & 0 & 0 & 0 & 0 & 0 & 0 & 1 & 0\\
0 & 0 & 0 & 0 & 1 & 1 & 1 & 0 & 0 & 1 & 1 & 0 & 0 & 1 & 1 & 0 & 0 & 0 & 1 & 1 & 1 & 0 & 1 & 0 & 1 & 0 & 1 & 0 & 1 & 0\\
0 & 0 & 0 & 0 & 0 & 0 & 0 & 0 & 1 & 0 & 0 & 0 & 1 & 0 & 0 & 1 & 0 & 0 & 1 & 1 & 0 & 0 & 0 & 0 & 0 & 0 & 1 & 0 & 1 & 0\\
0 & 0 & 0 & 0 & 0 & 1 & 1 & 0 & 1 & 0 & 0 & 0 & 1 & 0 & 0 & 0 & 0 & 1 & 0 & 0 & 0 & 1 & 1 & 0 & 0 & 0 & 1 & 0 & 0 & 0\\
0 & 1 & 1 & 0 & 0 & 0 & 0 & 0 & 0 & 0 & 0 & 0 & 1 & 0 & 1 & 0 & 1 & 1 & 1 & 0 & 1 & 0 & 1 & 1 & 1 & 1 & 0 & 0 & 0 & 0\\
0 & 0 & 0 & 0 & 1 & 1 & 1 & 0 & 0 & 1 & 1 & 1 & 0 & 0 & 0 & 0 & 0 & 0 & 1 & 0 & 0 & 0 & 0 & 0 & 0 & 1 & 1 & 0 & 0 & 1\\
0 & 0 & 0 & 0 & 1 & 1 & 0 & 1 & 1 & 0 & 0 & 0 & 0 & 0 & 0 & 0 & 1 & 1 & 0 & 1 & 1 & 0 & 1 & 1 & 1 & 0 & 0 & 0 & 1 & 0\\
0 & 1 & 0 & 1 & 1 & 0 & 1 & 0 & 1 & 0 & 0 & 1 & 0 & 0 & 0 & 1 & 0 & 1 & 1 & 1 & 0 & 1 & 1 & 0 & 1 & 0 & 0 & 0 & 1 & 1\\
0 & 1 & 0 & 0 & 0 & 0 & 0 & 0 & 0 & 1 & 0 & 0 & 0 & 0 & 1 & 0 & 0 & 0 & 0 & 1 & 1 & 0 & 1 & 0 & 1 & 1 & 0 & 1 & 1 & 1\\
0 & 0 & 1 & 0 & 1 & 0 & 1 & 1 & 0 & 0 & 0 & 1 & 0 & 1 & 0 & 0 & 0 & 1 & 0 & 1 & 1 & 0 & 1 & 0 & 1 & 0 & 1 & 1 & 1 & 1\\
1 & 1 & 1 & 0 & 1 & 0 & 0 & 0 & 0 & 0 & 1 & 1 & 0 & 1 & 1 & 0 & 1 & 0 & 0 & 0 & 1 & 1 & 0 & 1 & 1 & 1 & 1 & 0 & 0 & 0\\
1 & 0 & 0 & 1 & 0 & 0 & 1 & 0 & 1 & 1 & 0 & 1 & 1 & 0 & 1 & 0 & 0 & 0 & 0 & 1 & 0 & 1 & 1 & 1 & 0 & 0 & 0 & 1 & 1 & 0\\
1 & 0 & 0 & 0 & 0 & 1 & 0 & 0 & 1 & 1 & 0 & 0 & 0 & 1 & 1 & 1 & 1 & 0 & 1 & 0 & 1 & 0 & 1 & 1 & 0 & 1 & 0 & 1 & 0 & 1\\
1 & 1 & 0 & 1 & 0 & 0 & 0 & 0 & 1 & 0 & 0 & 1 & 0 & 1 & 0 & 1 & 1 & 1 & 0 & 1 & 0 & 1 & 1 & 1 & 1 & 0 & 0 & 1 & 1 & 0\\
1 & 0 & 1 & 0 & 1 & 0 & 0 & 0 & 0 & 0 & 1 & 0 & 0 & 0 & 1 & 0 & 0 & 1 & 1 & 0 & 1 & 0 & 1 & 0 & 1 & 0 & 1 & 0 & 1 & 1\\
1 & 1 & 0 & 0 & 1 & 1 & 0 & 0 & 1 & 0 & 1 & 1 & 0 & 1 & 1 & 1 & 1 & 0 & 1 & 1 & 1 & 1 & 0 & 0 & 1 & 1 & 1 & 1 & 1 & 1\\
0 & 1 & 1 & 0 & 1 & 0 & 1 & 0 & 0 & 0 & 0 & 1 & 0 & 1 & 0 & 0 & 0 & 1 & 1 & 1 & 1 & 0 & 0 & 0 & 1 & 0 & 0 & 1 & 0 & 1\\
0 & 1 & 0 & 0 & 0 & 1 & 0 & 0 & 1 & 0 & 0 & 1 & 0 & 1 & 1 & 1 & 1 & 1 & 0 & 0 & 1 & 1 & 1 & 1 & 0 & 0 & 1 & 1 & 1 & 1\\
0 & 0 & 1 & 0 & 0 & 0 & 1 & 0 & 0 & 0 & 0 & 1 & 1 & 0 & 0 & 1 & 0 & 1 & 0 & 1 & 0 & 0 & 1 & 0 & 0 & 0 & 1 & 1 & 1 & 0\\
1 & 0 & 0 & 0 & 0 & 0 & 0 & 0 & 1 & 1 & 1 & 0 & 1 & 0 & 0 & 0 & 1 & 1 & 0 & 0 & 0 & 1 & 1 & 0 & 1 & 1 & 0 & 1 & 1 & 1\\
1 & 0 & 0 & 0 & 1 & 0 & 0 & 0 & 0 & 0 & 0 & 0 & 0 & 0 & 0 & 1 & 1 & 0 & 1 & 1 & 1 & 0 & 1 & 1 & 1 & 1 & 1 & 0 & 1 & 1\\
1 & 0 & 1 & 1 & 0 & 0 & 1 & 1 & 1 & 1 & 0 & 0 & 0 & 1 & 1 & 1 & 1 & 0 & 1 & 0 & 1 & 1 & 1 & 0 & 1 & 1 & 1 & 1 & 0 & 1\\
0 & 0 & 0 & 0 & 1 & 1 & 1 & 0 & 0 & 0 & 0 & 0 & 1 & 0 & 1 & 1 & 1 & 0 & 0 & 1 & 0 & 1 & 1 & 1 & 1 & 0 & 1 & 1 & 1 & 0
\end{bmatrix}$
				\caption{An example of a graph $G$ where 
				$L_{\{k\}}(G)< \lfloor k \cdot Z^*_{rlx}(G)\rfloor$
				\label{ag}}
				\end{figure}

				In the sequel, we will prove that equality 
				$L_{\{k\}}(G)=k\cdot Z^*_{rlx}(G)$ holds for all graphs,
				but only for certain values of $k$.

				\begin{thm}\label{trmrel} For arbitrary graph $G$, 
					$(\exists q\in\mathbb{N})(\forall k_1\in\mathbb{N})\quad
					L_{\{k_1\cdot q\}}(G)=k_1\cdot q\cdot Z^*_{rlx}(G)$.
				\end{thm}
				\begin{pf}
					For arbitrary graph $G$, let $(x^*_1,\ldots,x^*_n)$ is an optimal solution 
					of linear programming formulation (\ref{LPrelobj})-(\ref{LPrelvar}), 
					with objective function value $Z^*_{rlx}(G)$.
					Since constraint matrix $A_G$
					is an integer matrix and right-hand side vector $b=(1\ 1 \ldots 1)^T$ is also the integer vector,
					then each feasible solution must be a vector with rational coordinates.
					Therefore, it also holds for optimal solution, i.e. $(\forall i)(x^*_i=\frac{p_i}{q_i}$ where
					$p_i\in \mathbb{Z}$, $q_i\in\mathbb{N}$ and $gcd(p_i,q_i)=1$ where 
					$gcd(a,b)$ is the greatest common divisor of $a$ and $b$. 
					Let us introduce $q=lcm(q_1,\ldots,q_n)$ where
					$lcm$ is the least common multiple. From the definition 
					it is obvious that $q_1,\ldots q_n\in \mathbb{N} \Rightarrow q\in \mathbb{N}$.
					If $x^*_i=0$ then $p_i=0$, let fix $q_i=1$ in that case.
					If (\ref{LPrelobj})-(\ref{LPrelvar}) has multiple optimal solutions we will 
					assume that we can arbitrarily choose one of them.\\

					Let $k=k_1\cdot q$ and let $(y^*_1,\ldots,y^*_n)$ is optimal solution
					of the dual problem of the linear programming formulation (\ref{LPrelobj})-(\ref{LPrelvar}). It 
					satisfies $A_G \cdot (y^*_1 \ldots y^*_n)^T \ge (1\ 1 \ldots 1)^T$.
					Since $(x^*_1,\ldots,x^*_n)$ and $(y^*_1,\ldots,y^*_n)$ are optimal solutions of the mutually dual problems it follows that values of corresponding objective functions are equal, that is $
					\sum\limits_{i = 1}^n {x_i^* }  = \sum\limits_{i = 1}^n {y_i^* } $. Dual problem of the problem (\ref{LPobj})-(\ref{LPvar}) is
					
					\begin{equation}
						\label{LDobj}
						\max \sum\limits_{i = 1}^{|V|} k \cdot {Y_i } = k \cdot \sum\limits_{i = 1}^{|V|} {Y_i }
					\end{equation}
					
					subject to

					\begin{equation}
						\sum\limits_{i = 1}^{|V|} {a_{ij} Y_i }  \ge 1,\quad j = 1, \ldots ,|V|
					\end{equation}
					
					\begin{equation}
						\label{LDvar}
						Y_i \ge 0  ,\quad i = 1, \ldots ,|V|
					\end{equation}
						
					As it can be seen value of objective function is equal to $k$ times of objective function of the dual of problem (\ref{LPrelobj})-(\ref{LPrelvar}). Now, it can be concluded that optimal value of objective function (\ref{LPobj}) is equal to $k \cdot \sum\limits_{i = 1}^n {x_i^* }$ and consequently that $(k \cdot x^*_1,\ldots,k \cdot x^*_n)$ is optimal solution 
					of linear programming formulation (\ref{LPobj})-(\ref{LPvar}). 
					As $k=q\cdot k_1$, such that $q=lcm(q_1,\ldots,q_n)$ and $(\forall i)x^*_i=\frac{p_i}{q_i}$
					then $k_1\cdot q\cdot x^*_i=k_1\cdot q\cdot \frac{p_i}{q_i}\in \mathbb{Z}$. 
					Since $(k_1\cdot q\cdot x^*_1,\ldots,k_1\cdot q\cdot x^*_n)$ is vector of integers,
					and it is optimal solution of linear programming formulation (\ref{LPobj})-(\ref{LPvar})
					then it is also optimal solution of integer linear programming formulation (\ref{ILPobj})-(\ref{ILPvar})
					with optimal value $k_1\cdot q \cdot Z^*_{rlx}$. 
					Therefore, $L_{\{k_1\cdot q\}}(G) = k_1\cdot q \cdot Z^*_{rlx}$ which confirms the statement of the theorem.
				\end{pf}

				\begin{cor}
					$\varlimsup\limits_{k\rightarrow +\infty} \frac{L_{\{k\}}(G)}{k}=Z^*_{rlx}(G)$
				\end{cor}
				\begin{pf}
					For a given graph $G$ let us consider sequence $(L_{\{k\}}(G))_{k\in \mathbb{N}}$ 
					and its subsequence $(L_{\{l\cdot q\}}(G))_{l\in \mathbb{N}}$ and
					$q\in\mathbb{N}$ as defined in Theorem \ref{trmrel}. 
					From Property \ref{prop1} follows that 
					$(\forall k)L_{\{k\}}(G) \le k\cdot Z^*_{rlx}(G)$ implying
					$(\forall k)\frac{L_{\{k\}}(G)}{k} \le Z^*_{rlx}(G)$.
					For subsequence $(L_{\{l\cdot q\}}(G))_{l\in \mathbb{N}}$
					from Theorem \ref{trmrel} it holds 
					$(\forall l)L_{\{l\cdot q\}}(G) = l\cdot q\cdot Z^*_{rlx}(G)$,
					so $(\forall l)\frac{L_{\{l\cdot q\}}(G)}{l\cdot q} = Z^*_{rlx}(G)$,
					implying 
					$\varlimsup\limits_{l\rightarrow +\infty} \frac{L_{\{l\cdot q \}}(G)}{l\cdot q}=Z^*_{rlx}(G)$,
					which directly confirms the statement.
				\end{pf}
				
				\begin{cor}
					For any graph $G$ there exists $q\in \mathbb{N}$ such that 
					$L_{\{k_1\cdot q\}}(G)$ can be found in polynomial time 
					for any $k_1\in \mathbb{N}$. 
				\end{cor}
				\begin{pf}
					Let us consider $q$ as defined in Theorem \ref{trmrel}. 
					If $k=q\cdot q_1$ then by Theorem \ref{trmrel}, optimal solution 
					of $L_{\{k\}}(G)$ can be obtained as optimal solution 
					of linear programming formulation 
					(\ref{LPobj})-(\ref{LPvar}).  
					Since it can be achieved in polynomial time,
					then in this case $L_{\{k\}}(G)$ can
					be obtained in polynomial time.
				\end{pf}
				
				\begin{obs}
				It should be noted that in Theorem \ref{np} (\cite{dob17}) word ''fixed'' is necessary. 
				Although for each simple connected graph $G$ and for some values of $k$, $L_{\{k\}(G)}$ can be
				determined in polynomial time, considered problem is still NP-complete for $k$ fixed.
				\end{obs}
				
				\begin{obs}
					It should be noted that $q$ defined in Theorem \ref{trmrel} is not necessarily minimal in the case
					with multiple optimal solution of (\ref{LPrelobj})-(\ref{LPrelvar}). The number of optimal solutions can
					be in worst case infinite (even uncountable), though all have the same optimal value, the minimal value
					of $q$ defined in Theorem \ref{trmrel} may not be obtained in polynomial time.\\
					Even in the case with single optimal solution of (\ref{LPrelobj})-(\ref{LPrelvar}), 
					$q=lcm(q_1,\cdots q_n)$ may not be the minimal $k$ for which 
					(\ref{LPrelobj})-(\ref{LPrelvar}) has integer optimal solution.
				\end{obs}
				
				Previous considerations were based on the Integer Linear Programming formulation of the proposed problem
				and its relaxation. Now, let us present several properties of $\{k\}$-packing function problem which are not derived from 
				ILP formulation. 
				In the following proposition, it will be proven that $\{1\}$-packing function problem of an arbitrary graph $G$
				can be reduced to vertex independence number problem on a graph $G^2$.

				\begin{prp}
					\label{thm:f1}
					$L_{\{1\}}(G)=\alpha(G^2)$.
				\end{prp}
				
				\begin{pf}
					($\Rightarrow$)
					Let $f$ be a $1$-packing function whose value $f(V(G))=L_{\{1\}}(G)$. 
					We define $I=\{v\in V(G)\,|\,f(v)=1\}$. 
					Let $u, v\in V(G)$, $u\ne v$ and $(u,v)\in E(G^2)$, i.e.
					$d(u,v)\le 2$.
					Then we have two cases:\\
					case 1:$v\in N(u)$. Since $f$ is $1$-packing function then 
					$f(N[u])=\sum\limits_{v\in  N[u])} f(v)\le 1$ implying
					$f(u)+f(v)\le 1$.\\
					case 2:$u,v\in N(w)$. Since 
					$f$ is $1$-packing function then 
					$f(N[w])=\sum\limits_{v\in  N[w])} f(v)\le 1$ implying
					$f(u)+f(v)\le 1$.\\
					In both cases we have $f(u)+f(v)\le 1$
					implying that $(u\notin I\lor v\notin I)$.
					Since for each edge from $E(G^2)$ has at least one endpoint in
					$I$, then $I$ is independent set of $G^2$.
					
					($\Leftarrow$)
					Let $I$ be an independent set of $G^2$. We define $f(v) = \begin{cases}
					1, \,v \in I\\
					0, \,v \notin I\\
					\end{cases}. $

					Let $v$ be an arbitrary vertex from $V(G)$, 
					and $u,w\in N(v)$ and 
					$u\ne w$. Then, $d(u,w)\le 2$. Since $I$ is an independent set of 
					$G^2$ at most one of vertices $u, w$ is in $I$, so $f(u)+f(v)+f(w)\le 1$. 
					Since $u$ and $w$ are arbitrary vertices from $N(v)$, 
					then $f(N[v])=\sum\limits_{w\in  N[v])} f(w)\le 1$. 
					In the case when $v$ has only one neighbor $u$,
					holds $f(N[v])=f(u)+f(v)\le 1$.
					Since $v$ is an arbitrary vertex from $V(G)$ it follows that
					$f$ is $1$-packing function of $G$.

				\end{pf}

				\begin{cor}
					$L_{\{1\}}(G)=\rho_2(G)$
				\end{cor}
				
				\begin{cor} If $Diam(G)=2$, then $L_{\{1\}}(G)=1$. 
				\end{cor}
				\begin{pf}
					If $Diam(G)=2$, then $G^2=K_{|V(G)|}$, and consequently, $L_{\{1\}}(G)= \alpha(K_{|V(G)|})=1$.
				\end{pf}

				Next, it will be proposed computationally simple lower bound based upon the graph diameter.
				
				\begin{prp}
					$L_{\{k\}}(G)\ge \lceil\frac{1+Diam(G)}{3}\rceil\cdot k$
				\end{prp}
				
\begin{pf}
				From Proposition \ref{m1} it stands that 
				$L_1(G)\ge \lceil\frac{diam(G)+1}{3}\rceil$. 
				On the other hand, from Proposition \ref{a1} it stands that $L_{\{k\}}\ge k\cdot L_1$. 
				By combining mentioned inequalities we obtain
				$L_{\{k\}}\ge k\cdot L_1 \ge k\cdot \lceil\frac{diam(G)+1}{3}\rceil$ 
\end{pf}
				
				This lower bound is tight as it can be seen from Proposition \ref{path}.\\

				Next, it will be introduced upper bound based on the vertices' degree.
				
				\begin{prp}
					\label{prp1}
					$L_{\{k\}}(G) \le \lfloor \frac{nk}{1+\delta(G)}\rfloor$.
				\end{prp}
				\begin{pf}
					For each vertex $v\in V(G)$ it holds that
					$f(N[v])\le k$. Summing previous inequalities on all
					vertices from $V$ we obtain:
					$$n\cdot k\ge \sum\limits_{v\in V}f(N[v])=\sum\limits_{v\in V}\sum\limits_{w\in N[v]}f(w)$$.
					  
					On the other hand, for arbitrary vertex $u$ from $V$,
					in previous sums $f(u)$ appears exactly $1+deg(u)$ times:
					once for vertex $u$ and $deg(u)$ times for each vertex that is adjacent to the vertex $u$. Therefore, we get:
					$$\sum\limits_{v\in V}\sum\limits_{w\in N[v]}f(w)=
					\sum\limits_{u\in V}(1+deg(u))\cdot f(u)\ge \sum\limits_{u\in V}(1+\delta)\cdot f(u)=$$
					$$=(1+\delta)\cdot\sum\limits_{u\in V} f(u)=(1+\delta)\cdot f(V(G))$$.  					
					As a consequence, it holds 
					$$f(V(G))\le \frac{n\cdot k}{1+\delta}\Rightarrow 
					f(V(G))\le \left\lfloor\frac{n\cdot k}{1+\delta}\right\rfloor$$ 
					The previous inequality holds because $f(V(G))\in \mathbb{N}\cup\{0\}$.
				\end{pf}
				
				\begin{cor}
					\label{regular}
					If $G$ is a regular graph of degree $r$,  then $L_{\{k\}}(G) \le  \lfloor \frac{nk}{1+r}\rfloor$
				\end{cor}
				
				Bounds in Proposition \ref{prp1} are tight as it can be seen from the two following statements.

				\begin{prt}
					For complete graph (clique) $K_n$ holds $L_{\{k\}}(K_n)=k$.
				\end{prt}

				\begin{prp} 
					For cycle $C_n$ holds
					$L_{\{k\}}(C_n)= \lfloor \frac{n\cdot k}{3}\rfloor$.
				\end{prp}
				\begin{pf}
					Let graph $C_n$ be a cycle, i.e. 
					$C_n=(V,E)$ where $V=\{0,1,2,\ldots,n-1\}$ and
					$E=\{\{0,1\},\{1,2\},\{2,3\},\ldots,\{n-2,n-1\},\{n-1,0\}\}$.
					
					Let us define function $f$ as
					
					$$
					f(v_i)= \left\{
					\begin{array}{ll}
						\lfloor\frac{k}{3}\rfloor,     & i\equiv 0\pmod 3, \\
						\lfloor\frac{k}{3}+0.5\rfloor, & i\equiv 1\pmod 3, \\
						\lceil\frac{k}{3}\rceil,       & i\equiv 2\pmod 3. 
					\end{array}
					\right.
					$$
					All possible cases are presented in Table \ref{cn}
					
					\begin{table}
					{\setlength\tablinesep{4pt}
						\tiny
					    \centering
						\caption{$f(N[v])$ for $C_n$}\label{cn}
					\begin{tabular}{|c|c|c|c|}
						
						\hline 
						$n$    & $k$    & $v$ & $f(N[v])$ \\ 
						\hline 
						      $3m$ & $3l$ & $v_i, i=0,\ldots,3m-1$ & $\lfloor\frac{3l}{3}\rfloor+\lfloor\frac{3l}{3}+0.5\rfloor+ \lceil\frac{3l}{3}\rceil=l+l+l=3l=k\le k$ \\				       
						\hline
						\multicolumn{4}{|c|}{$f(V(C_{3m}))=m\cdot3l=\lfloor\frac{nk}{3}\rfloor$}\\ 
						\hline
						
						$3m$   & $3l+1$ & $v_i, i=0,\ldots,3m-1$ & 
						$\lfloor\frac{3l+1}{3}\rfloor+\lfloor\frac{3l+1}{3}+0.5\rfloor+ \lceil\frac{3l+1}{3}\rceil=l+l+l+1=3l+1=k\le k$\\
						\hline					
						\multicolumn{4}{|c|}{$f(V(C_{3m}))=m\cdot(3l+1)=\lfloor\frac{nk}{3}\rfloor$}\\
						\hline	

						$3m$   & $3l+2$ & $v_i, i=0,\ldots,3m-1$ & 
						$\lfloor\frac{3l+2}{3}\rfloor+\lfloor\frac{3l+2}{3}+0.5\rfloor+ \lceil\frac{3l+2}{3}\rceil=l+l+1+l+1=3l+2=k\le k$\\
						\hline			
						\multicolumn{4}{|c|}{$f(V(C_{3m}))=m\cdot(3l+2)=\lfloor\frac{nk}{3}\rfloor$}\\
						\hline
						\hline

						\hline											
						$3m+1$ & $3l$ & $v_0$ &     $\lfloor\frac{3l}{3}\rfloor+\lfloor\frac{3l}{3}\rfloor+\lfloor\frac{3l}{3}+0.5\rfloor=l+l+l=3l=k\le k$\\
						\hline
						$3m+1$ & $3l$ & $v_i, i=1,\ldots,3m-1$ &     $\lfloor\frac{3l}{3}\rfloor+\lfloor\frac{3l}{3}+0.5\rfloor+\lceil\frac{3l}{3}\rceil=l+l+l=3l=k\le k$\\
						\hline
						$3m+1$ & $3l$ & $v_{3m}$ &     $\lceil\frac{3l}{3}\rceil+\lfloor\frac{3l}{3}\rfloor+\lfloor\frac{3l}{3}\rfloor=l+l+l=3l=k\le k$\\
						\hline
						\multicolumn{4}{|c|}{$f(V(C_{3m+1}))=(m+1)\cdot l+ml+ml=3ml+l=l(3m+1)=\lfloor\frac{(3m+1)3l}{3}\rfloor=\lfloor\frac{nk}{3}\rfloor$}\\							\hline
						\hline
																								
						$3m+1$ & $3l+1$ & $v_0$ & 
					    $\lfloor\frac{3l+1}{3}\rfloor+\lfloor\frac{3l+1}{3}\rfloor+\lfloor\frac{3l+1}{3}+0.5\rfloor=l+l+l=3l=k-1\le k$\\
						\hline
						$3m+1$ & $3l+1$ & $v_i, i=1,\ldots,3m-1$ &     $\lfloor\frac{3l+1}{3}\rfloor+\lfloor\frac{3l+1}{3}+0.5\rfloor+\lceil\frac{3l+1}{3}\rceil=l+l+l+1=3l+1=k\le k$\\																																										\hline
						$3m+1$ & $3l+1$ & $v_{3m}$ &     $\lceil\frac{3l+1}{3}\rceil+\lfloor\frac{3l+1}{3}\rfloor+\lfloor\frac{3l+1}{3}\rfloor=l+1+l+l=3l+1=k\le k$\\
						\hline
						\multicolumn{4}{|c|}{$f(V(C_{3m+1}))=(m+1)\cdot l+ml+m(l+1)=3ml+m+l=\lfloor\frac{(3m+1)(3l+1)}{3}\rfloor=\lfloor\frac{nk}{3}\rfloor$}\\
						\hline
						\hline

						$3m+1$ & $3l+2$ & $v_0$ & $\lfloor\frac{3l+2}{3}\rfloor+\lfloor\frac{3l+2}{3}\rfloor+\lfloor\frac{3l+2}{3}+0.5\rfloor=l+l+l+1=3l+1=k-1\le k$\\
						\hline
						$3m+1$ & $3l+2$ & $v_i, i=1,\ldots,3m-1$ &     $\lfloor\frac{3l+2}{3}\rfloor+\lfloor\frac{3l+2}{3}+0.5\rfloor+\lceil\frac{3l+2}{3}\rceil=l+l+1+l+1=3l+2=k\le k$\\																									\hline
						$3m+1$ & $3l+2$ & $v_{3m}$ &     $\lceil\frac{3l+2}{3}\rceil+\lfloor\frac{3l+2}{3}\rfloor+\lfloor\frac{3l+2}{3}\rfloor=l+1+l+l=3l+1=k-1\le k$\\
						\hline
						\multicolumn{4}{|c|}{$f(V(C_{3m+1}))=(m+1)\cdot l+m(l+1)+m(l+1)=3ml+l+2m=\lfloor\frac{(3m+1)(3l+2)}{3}\rfloor=\lfloor\frac{nk}{3}\rfloor$}\\
\hline						
						\hline	
						$3m+2$ & $3l$ & $v_0$ &     $\lfloor\frac{3l}{3}+0.5\rfloor+\lfloor\frac{3l}{3}\rfloor+\lfloor\frac{3l}{3}+0.5\rfloor=l+l+l=3l=k\le k$\\
						\hline
						$3m+2$ & $3l$ & $v_i, i=1,\ldots,3m$ &     $\lfloor\frac{3l}{3}\rfloor+\lfloor\frac{3l}{3}+0.5\rfloor+\lceil\frac{3l}{3}\rceil=l+l+l=3l=k\le k$\\
						\hline
						$3m+2$ & $3l$ & $v_{3m+1}$ &     $\lfloor\frac{3l}{3}\rfloor+\lfloor\frac{3l}{3}+0.5\rfloor+\lfloor\frac{3l}{3}\rfloor=l+l+l=3l=k\le k$\\
						\hline
						\multicolumn{4}{|c|}{$f(V(C_{3m+2}))=(m+1)\cdot l+(m+1)\cdot l+m\cdot l=3ml+2l=\lfloor\frac{(3m+2)3l}{3}\rfloor=\lfloor\frac{nk}{3}\rfloor$}\\
						\hline
						\hline
						$3m+2$ & $3l+1$ & $v_0$ & $\lfloor\frac{3l+1}{3}+0.5\rfloor+\lfloor\frac{3l+1}{3}\rfloor+\lfloor\frac{3l+1}{3}+0.5\rfloor=l+l+l=3l=k-1\le k$\\
						\hline
						$3m+2$ & $3l+1$ & $v_i, i=1,\ldots,3m$ &     $\lfloor\frac{3l+1}{3}\rfloor+\lfloor\frac{3l+1}{3}+0.5\rfloor+\lceil\frac{3l+1}{3}\rceil=l+l+l+1=3l+1=k\le k$\\																									\hline
						$3m+2$ & $3l+1$ & $v_{3m+1}$ &     $\lfloor\frac{3l+1}{3}\rfloor+\lfloor\frac{3l+1}{3}+0.5\rfloor+\lfloor\frac{3l+1}{3}\rfloor=l+l+l=3l=k-1\le k$\\
						\hline
						\multicolumn{4}{|c|}{$f(V(C_{3m+2}))=(m+1)\cdot l+(m+1)\cdot l+m(l+1) = 3ml+m+2l=\lfloor\frac{(3m+2)(3l+1)}{3}\rfloor= \lfloor\frac{nk}{3}\rfloor$}\\						\hline
						\hline						

						$3m+2$ & $3l+2$ & $v_0$ & $\lfloor\frac{3l+2}{3}+0.5\rfloor+\lfloor\frac{3l+2}{3}\rfloor+\lfloor\frac{3l+2}{3}+0.5\rfloor=l+1+l+l+1=3l+2=k\le k$\\
						\hline
						$3m+2$ & $3l+2$ & $v_i, i=1,\ldots,3m$ &     $\lfloor\frac{3l+2}{3}\rfloor+\lfloor\frac{3l+2}{3}+0.5\rfloor+\lceil\frac{3l+2}{3}\rceil=l+l+1+l+1=3l+2=k\le k$\\																									\hline
						$3m+2$ & $3l+2$ & $v_{3m+1}$ &     $\lfloor\frac{3l+2}{3}\rfloor+\lfloor\frac{3l+2}{3}+0.5\rfloor+\lfloor\frac{3l+2}{3}\rfloor=l+l+1+l=3l+1=k-1\le k$\\
						\hline \multicolumn{4}{|c|}{$f(V(C_{3m+2}))=
							(m+1)\cdot l + (m+1)(l+1) + m(l+1) = 3ml+2l+2m+1=\lfloor\frac{(3m+2)(3l+2)}{3}\rfloor =\lfloor\frac{nk}{3}\rfloor$}\\																\hline

\hline 						
					\end{tabular} }
					\end{table}

					From Table it is obvious that in each case $f(N[w])\le k$ and $f(V(G))=\lfloor\frac{n\cdot k}{3}\rfloor$. Therefore we proved that $L_{\{k\}}(C_n)\ge \lfloor\frac{n\cdot k}{3}\rfloor$.
					Since $C_n$ is regular graph with $r=2$ it holds that 
					$L_{\{k\}}(G) \le  \lfloor \frac{nk}{1+2}\rfloor=\lfloor \frac{nk}{3}\rfloor$.
					Consequently, equality 
					$L_{\{k\}}(G) = \lfloor \frac{nk}{3}\rfloor$ holds.
				\end{pf}
				
				\section{Conclusions}
				In this paper the $\{k\}$-packing function problem is
				studied. First, special relation  was established between cases when $k=1$, $k\ge 2$, 
				and the optimal solution of the linear programming relaxation. Second, sufficient conditions for optimality were introduced. 
				It was proven that, for arbitrary simple connected graph $G$ and some values of $k$,
				$L_{\{k\}}(G)$ can be determined in the polynomial time.
				Next, $\{1\}$-packing function problem was studied and its connection
				with the independent set number and $2$-packing problem. 
				Finally, lower and upper bound was introduced as well as 
				optimal values for some special classes of graphs. 
				
				The future work could be directed to considering the $\{k\}$-packing function number 
				of some challenging classes of graphs.

				\subsection*{Acknowledgements}
				This research was partially supported by Serbian Ministry of
				Education, Science and Technological Development under the grants no. 174010 and
				174033.
				

\end{document}